\ifdefined\PreambleLoaded\else
  \documentclass{amsart}
\usepackage{marktext}
\usepackage{gimac}
\RequirePackage[nointlimits]{esint}

\newcommand{\eff}{\mathrm{eff}}

\newcommand{\std}{\mathrm{std}}

\newcommand{\ttmix}{\tilde{t}_\mathrm{mix}}
\newcommand{\e}{\bm{e}}

\DeclareMathOperator{\dist}{dist}

\DeclareMathOperator{\cov}{cov}
\DeclareMathOperator{\unif}{unif}

\DeclarePairedDelimiter{\ip}{\langle}{\rangle}

\newtheorem{assumption}[theorem]{Assumption}

\colorlet{darkgreen}{green!50!black}

\newcommand{\PreambleLoaded}{}

\fi
\RequirePackage{tikz}
\begin{document}
\title{Residual Diffusivity for Expanding Bernoulli Maps}
\begin{abstract}
  Consider a discrete time Markov process~$X^\epsilon$ on~$\R^d$ that makes a deterministic jump based on its current location, and then takes a small Gaussian step of variance~$\epsilon^2$.
  We study the behavior of the \emph{asymptotic variance} as~$\epsilon \to 0$.
  In some situations (for instance if there were no jumps), then the asymptotic variance vanishes as~$\epsilon \to 0$.
  When the jumps are ``chaotic'', however, the asymptotic variance may be bounded from above and bounded away from~$0$, as~$\epsilon \to 0$.
  This phenomenon is known as \emph{residual diffusivity}, and we prove this occurs when the jumps are determined by certain expanding Bernoulli maps.
\end{abstract}
\author[Cooperman]{William Cooperman}
\address{%
  Department of Mathematics, ETH Z\"urich, Switzerland
}
\email{bill@cprmn.org}
\author[Iyer]{Gautam Iyer}
\address{%
  Department of Mathematical Sciences, Carnegie Mellon University, Pittsburgh, PA 15213.
}
\email{gautam@math.cmu.edu}
\author[Nolen]{James Nolen}
\address{%
  Department of Mathematics, Duke University, Durham, NC 27708.
}
\email{james.nolen@duke.edu}
\thanks{This work has been partially supported by the National Science Foundation under grants
  DMS-2303355, 
  DMS-2342349, 
  DMS-2406853, 
  and the Center for Nonlinear Analysis.}

\maketitle

\section{Introduction}\label{s:intro}

\subsection{Main Result}
Consider a Markov process $\{X^\epsilon_n\}_{n\geq 0}$ that makes a deterministic jump based on its current location, followed by a Gaussian step of variance~$\epsilon^2$.
Explicitly, $X^\epsilon_{n+1}$ is determined from~$X^\epsilon_n$ by
\begin{equation}\label{e:Xn}
  X^\epsilon_{n+1} = \varphi(X^\epsilon_n) + \epsilon \xi_{n+1} \,.
\end{equation}
Here $\{\xi_n\}_{n\geq 1}$ is a family of independent standard Gaussian random variables, and~$\varphi \colon \R^d \to \R^d$ is a Lebesgue measure preserving map with a periodic displacement (i.e.\ the function~$x \mapsto \varphi(x) - x$ is~$\Z^d$ periodic).

In this situation one can see that for any~$v \in \R^d$, the variance~$\var( v \cdot X^\epsilon_n )$ grows linearly with $n$ as~$n \to \infty$.
Our interest is to study the growth rate of the variance asymptotically as~$\epsilon \to 0$.
More precisely, we are interested in the behavior of the \emph{asymptotic variance}
\begin{equation}\label{e:AV}
    \lim_{n \to \infty} \frac{\var^{\mu_0} (v \cdot X^\epsilon_n)}{n}
\end{equation}
in the vanishing noise limit~$\epsilon \to 0$.
Here~$\mu_0$ is a probability distribution on~$\R^d$, and the notation~$\var^{\mu_0}(v \cdot X^\epsilon_n)$ denotes the variance of~$v \cdot X^\epsilon_n$ given~$X^\epsilon_0 \sim \mu_0$.

If~$\varphi$ is ``not too chaotic'', then the asymptotic variance of~$v \cdot X^\epsilon$ may either vanish as~$\epsilon \to 0$ (for instance, if~$\varphi$ is the identity map) or diverge to $+\infty$ (for instance, if the map $\varphi$ has a shear structure, with unbounded orbits).
If, however, $\varphi$ is ``chaotic'', then it may be possible for the asymptotic variance of~$v \cdot X^\epsilon$ to be bounded from above and bounded away from~$0$ as~$\epsilon \to 0$.
This phenomenon is known as \emph{residual diffusivity} and its study originated in fluid dynamics~\cite{Taylor21,BiferaleCrisantiEA95,MurphyCherkaevEA17}.

The purpose of this paper is to produce a deterministic example of~$\varphi$ for which~$X^\epsilon$ exhibits residual diffusivity.
The difficulty in rigorously proving residual diffusivity is that it involves two non-interchangeable limits: the vanishing noise limit, and the long time limit.
A natural candidate for maps that exhibit residual diffusivity are ones where the interchanged limit (i.e.\ the asymptotic variance of~\eqref{e:Xn} with~$\epsilon = 0$) is non-zero, at least for some initial distribution (say a small bump function).
This is easy to compute for several examples such as the doubling map, toral automorphisms, bakers map, and expanding Bernoulli maps.
Of these, the only one we are able to prove residual diffusivity for is a class of expanding Bernoulli maps.
This is the main result of this paper.
To the best of our knowledge, this is the only class of chaotic maps for which residual diffusion has been rigorously proved.
\begin{figure}[hbt]
  \resizebox{.8\linewidth}{!}{
    \begin{tikzpicture}
      \draw[gray, thick] (0, 0) -- (0, 6) -- (6, 6) -- (6, 0) -- (0, 0);
      \draw[gray, thick] (6.5, 0) -- (6.5, 6) -- (12.5, 6) -- (12.5, 0) -- (6.5, 0);
      \draw[gray] (2, 0) -- (2, 6);
      \draw[gray] (4, 0) -- (4, 6);
      \draw[gray] (8.5, 0) -- (8.5, 6);
      \draw[gray] (10.5, 0) -- (10.5, 6);
      \draw[gray] (0, 2) -- (6, 2);
      \draw[gray] (0, 4) -- (6, 4);
      \draw[gray] (6.5, 2) -- (12.5, 2);
      \draw[gray] (6.5, 4) -- (12.5, 4);

      \draw (2, 3) .. controls (2.3, 2.3) .. (3, 2.5) .. controls (3.5, 2.2) .. (4, 3);
      \draw (3, 4) .. controls (3, 3) .. (4, 3.5);
      \draw (4, 2.3) .. controls (3, 2.1) and (3, 2) .. (3, 2.5) .. controls (3.3, 2.6) and (3, 3) .. (2, 4);

      \path[fill=red] (2, 2) -- (4, 2) -- (4, 2.3) .. controls (3, 2.1) and (3, 2) .. (3, 2.5) .. controls (2.3, 2.3) .. (2, 3) -- (2, 2);
      \path[fill=green] (4, 2.3) .. controls (3, 2.1) and (3, 2) .. (3, 2.5) .. controls (3.5, 2.2) .. (4, 3);
      \path[fill=blue] (4, 3) .. controls (3.5, 2.2) .. (3, 2.5) .. controls (3.3, 2.6) and (3, 3) .. (2, 4) -- (3, 4) .. controls (3, 3) .. (4, 3.5) -- (4, 3);
      \path[fill=orange] (2, 3) .. controls (2.3, 2.3) .. (3, 2.5) .. controls (3.3, 2.6) and (3, 3) .. (2, 4) -- (2, 3);
      \path[fill=purple] (3, 4) .. controls (3, 3) .. (4, 3.5) -- (4, 4) -- (3, 4);

      \path[fill=red] (10.5, 4) -- (12.5, 4) -- (12.5, 6) -- (10.5, 6) -- (10.5, 4);
      \path[fill=green] (8.5, 2) -- (10.5, 2) -- (10.5, 4) -- (8.5, 4) -- (8.5, 2);
      \path[fill=blue] (10.5, 2) -- (12.5, 2) -- (12.5, 4) -- (10.5, 4) -- (10.5, 2);
      \path[fill=orange] (6.5, 0) -- (8.5, 0) -- (8.5, 2) -- (6.5, 2) -- (6.5, 0);
      \path[fill=purple] (6.5, 4) -- (8.5, 4) -- (8.5, 6) -- (6.5, 6) -- (6.5, 4);
    \end{tikzpicture}
  }
  \caption{
    One example of an expanding Bernoulli map~$\varphi$.
    The colored regions on the left are mapped to regions of the same color on the right.
  }
  \label{f:BernoulliExample}
\end{figure}
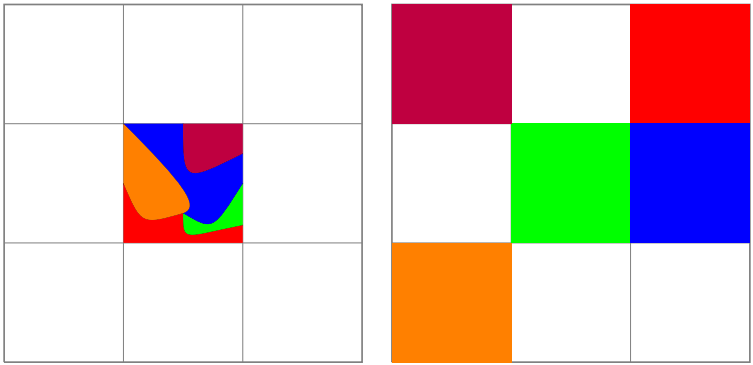

Our main result is the following:
\begin{theorem}\label{t:residual-diffusivity}
  Suppose~$\varphi$ is obtained from an expanding Bernoulli map satisfying the conditions in Assumption~\ref{a:phiTilde}, below.
  For all $v \in \R^d$, and all subgaussian initial distributions~$\mu_0$ we have
  %
  \begin{equation}\label{e:residual-diffusivity}
    \lim_{\epsilon \to 0}
      \lim_{n \to \infty} \frac{\var^{\mu_0}( v \cdot X^\epsilon_n )}{n}
	= \var^{\pi_0} (v \cdot \floor{X^0_1})
      \,.
  \end{equation}
\end{theorem}

Here $\pi_0$ denotes the uniform distribution on the unit cube~$Q_0 \defeq [0, 1)^d$, and the notation~$\floor{\cdot}$ denotes the bottom left vertex of the containing unit integer lattice cube.
Explicitly, for~$x \in \R^d$, the notation $\floor{x}$ denotes the unique element in~$\Z^d$ for which~$x \in \floor{x} + Q_0$.
We also recall that a distribution is subgaussian if there exists~$\lambda > 0$ such that~$\int_{\R^d} e^{\lambda |x|^2} d\mu_0(x) < \infty$.

\begin{figure}[hbt]
  \includegraphics[width=.5\linewidth]{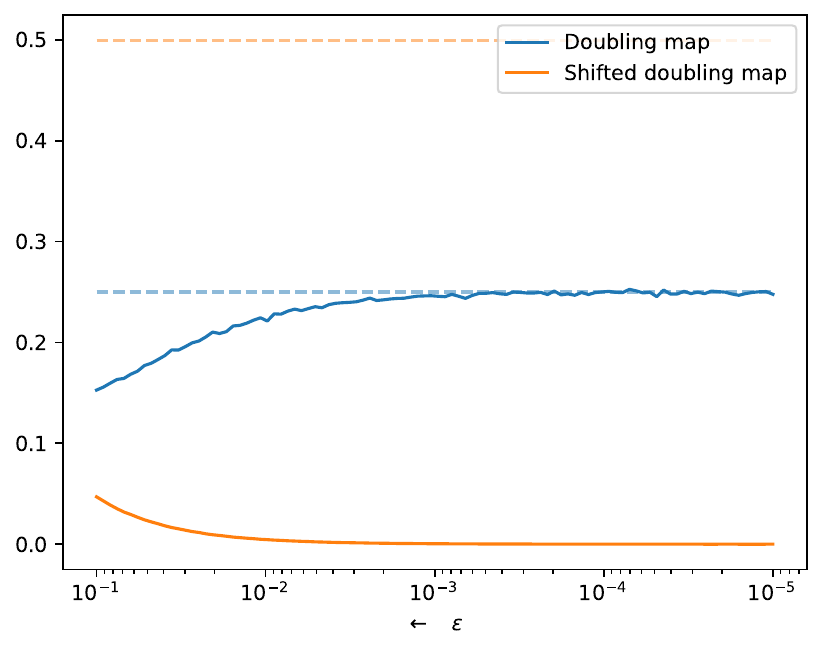}
  \caption{
    Asymptotic variance for the doubling map, and the shifted doubling map as~$\epsilon$ varies from~$10^{-1}$ to~$10^{-5}$.
    For the doubling map the asymptotic variance approaches the value on the right of~\eqref{e:residual-diffusivity} (blue dashed line) as~$\epsilon \to 0$.
    The shifted doubling map, however, doesn't exhibit residual diffusivity and the asymptotic variance approaches~$0$ instead of the value on the right of~\eqref{e:residual-diffusivity} (orange dashed line).
  }
  \label{f:1dNumerics}
\end{figure}

Notice that the right side of \eqref{e:residual-diffusivity} is non-zero if the image of $Q_0$ under $\varphi$ intersects more than one cube. One example of an expanding Bernoulli map that satisfies the conditions of Theorem~\ref{t:residual-diffusivity} is shown pictorially in Figure~\ref{f:BernoulliExample}.
Another example can be constructed from the one dimensional doubling map, and is defined by 	
\begin{equation}\label{e:doubling}
  \varphi(x) = 2x - \floor{x}
  \,,\quad x \in \R
  \,.
\end{equation} 
Interestingly, if we introduce a shift by~$1/2$ into~\eqref{e:doubling} then we obtain an example of a mixing map which does \emph{not} exhibit residual diffusivity.
Explicitly, we define the \emph{shifted doubling map} by
\begin{equation}\label{e:doubling2}
  \varphi_2(x) =
      2x - \floor{x} - \frac{1}{2}
    \,,\quad x \in \R
  \,.
\end{equation}
The shifted doubling map satisfies all but one of the required hypothesis in Theorem~\ref{t:residual-diffusivity} (specifically, for the shifted doubling map the cube centers~$o_i$ in Section~\ref{s:bernoulli}, below, are not integers).
The results of the numerical simulations in Figure~\ref{f:1dNumerics} show that the doubling map exhibits residual diffusivity, while the shifted doubling map does not.

The main idea behind the proof of Theorem~\ref{t:residual-diffusivity} is to use mixing to show that the asymptotic variance starting from any (subgaussian) initial distribution~$\mu_0$ is close to the asymptotic variance starting from the uniform distribution~$\pi_0$ on $Q_0$.
When~$X^\epsilon_0 \sim \pi_0$, the Bernoulli structure of~$\varphi$ serendipitously makes increments of~$\floor{\varphi(X^\epsilon_n)} - \floor{X^\epsilon_n}$ independent, allowing us to compute the asymptotic variance explicitly.
Before delving into the technical details of the proof we briefly survey the literature and place Theorem~\ref{t:residual-diffusivity} in the context of existing results.

\subsection{Motivation and Literature review}
Our interest in this problem stems from understanding the long time behavior of diffusions whose drift has ``chaotic trajectories''.
Explicitly, consider the continuous time diffusion process $X^\epsilon_t$ defined by the SDE
\begin{equation}\label{e:SDEX}
  d X^\epsilon_t = u(X^\epsilon_t) \, dt + \epsilon \, dW_t
  \quad\text{on } \R^d
  \,,
\end{equation}
where~$W$ is a~$d$-dimensional Brownian motion, and~$u$ is a spatially periodic, divergence free vector field.
One physical situation where this is relevant is in the study of diffusive tracer particles being advected by an incompressible fluid.

On small (i.e.\ $O(1)$) time scales, it the process of~$X^\epsilon$ stays close to the deterministic trajectories of~$u$, and a large deviations principle can be established (see for instance~\cite{FreidlinWentzell12}).
On intermediate (i.e.\ $O(\abs{\ln \epsilon}/\epsilon^\alpha)$ for~$\alpha \in [0, 2)$) time scales certain non-Markovian effects arise and lead to anomalous diffusion~\cite{
  Young88,
  YoungPumirEA89,
  Bakhtin11,
  HairerKoralovEA16,
  HairerIyerEA18
}.
On long (i.e.\ $O(1/\epsilon^2)$) time scales, homogenization occurs and the net effect of the drift can be averaged and the process~$X^\epsilon$ can be approximated by a Brownian motion with covariance matrix~${D^\epsilon_\eff}$ called the \emph{effective diffusivity}.
This was first studied in this setting by Freidlin~\cite{Freidlin64}, and is now the subject of many standard books with several important applications~\cite{BensoussanLionsEA78,PavliotisStuart08}.

The effective diffusivity matrix ${D^\epsilon_\eff}$ is the unique symmetric matrix whose action on vectors~$v \in \R^d$ is given by
\begin{equation}\label{e:Deff}
  v {D^\epsilon_\eff} v
     = \lim_{t \to \infty} \frac{ \var( v \cdot X^\epsilon_t ) }{t}
     \,.
\end{equation}
This, however, is hard to compute explicitly and authors usually characterize it in terms of a cell problem on $\T^d$, or a variational problem.
In a few special situations (such as shear flows, or cellular flows) the asymptotic behavior of~${D^\epsilon_\eff}$ as~$\epsilon \to 0$ is known~\cite{
  Taylor53,
  ChildressSoward89,
  FannjiangPapanicolaou94,
FannjiangPapanicolaou97,
MajdaKramer99,
  Heinze03,
  Koralov04,
NovikovPapanicolaouEA05,
  RyzhikZlatos07
}.

The motivation for the present paper comes from thinking about the case when the deterministic flow of~$u$ is chaotic on the torus.
In this case it has been conjectured that~${D^\epsilon_\eff}$ is~$O(1)$ as~$\epsilon \to 0$, a phenomenon known as \emph{residual diffusivity}.
Study of this was initiated by Taylor~\cite{Taylor21} over 100 years ago, and has since been extensively studied by many authors~\cite{
  ZaslavskyStevensEA93,
  BiferaleCrisantiEA95,
  MajdaKramer99,
  Zaslavsky02,
  MurphyCherkaevEA17
}.
While this has been confirmed numerically and studied for elephant random walks~\cite{
  LyuXinEA17,
  LyuXinEA18,
  MurphyCherkaevEA20,
  WangXinEA21,
  WangXinEA22,
  LyuWangEA22,
  KaoLiuEA22
},
we are unaware of a rigorous proof for any any example.
\medskip

The goal of this paper is to rigorously exhibit residual diffusivity in a simple setting.
First we replace the notion of ``chaotic'' with the assumption that the flow of~$u$ on the torus~$\T^d$ is \emph{exponentially mixing} (see~\cite{SturmanOttinoEA06}).
In continuous time, however, examples of exponentially mixing flows are not easy to construct.
The canonical example of an exponentially mixing flow is the geodesic flow on the unit sphere bundle of negatively curved manifolds~\cite{Dolgopyat98}. On the~$3$-torus, however, the existence of a divergence free, $C^1$, time independent, exponentially mixing velocity field is an open question.
To the best of our knowledge, there are only examples of lower regularity~\cite{ElgindiZlatos19}, and several time dependent examples~\cite{Pierrehumbert94,BedrossianBlumenthalEA19,MyersHillSturmanEA22,BlumenthalCotiZelatiEA23,ElgindiLissEA23,ChristieFengEA23}.

On the other hand, there are several simple, well known, examples of exponentially mixing dynamical systems on the torus~\cite{SturmanOttinoEA06}.
Therefore, instead of studying a continuous-time system, we study the discrete time system~\eqref{e:Xn}.
The system~\eqref{e:Xn} can be viewed as a discretization of~\eqref{e:SDEX} where~$\varphi$ is the flow map of~$u$ after a fixed amount of time.
The periodicity of~$u$ translates to the requirement that the displacement~$x \mapsto \varphi(x) - x$ is periodic.
Incompressibility of~$u$ translates to the requirement that~$\varphi$ is Lebesgue measure preserving.
This leads us to study~\eqref{e:Xn} in the general situation that~$\varphi$ is any Lebesgue measure preserving map with a periodic displacement (and not necessarily a diffeomorphism obtained as the flow of an incompressible vector field).
Such systems are interesting in their own right, and various aspects of them have been extensively studied~\cite{
  FannjiangWoowski03,
  FannjiangNonnenmacherEA04,
  ThiffeaultChildress03,
  FengIyer19,
  OakleyThiffeaultEA21,
  IyerLuEA24
}.

In this time-discrete setting, Theorem~\ref{t:residual-diffusivity} exactly computes the vanishing noise limit of the effective diffusivity.
Our proof, however, relies on the Bernoulli structure of~$\varphi$ and will not apply to general mixing maps.
For general mixing maps~$\varphi$ we conjecture that~\eqref{e:residual-diffusivity} in Theorem~\ref{t:residual-diffusivity} should be replaced with
\begin{equation}\label{e:residual-diffusivity2}
  \lim_{\epsilon \to 0}
  \lim_{n \to \infty} \frac{\var^{\mu_0}( v \cdot X^\epsilon_n )}{n}
  = \lim_{n \to \infty} \frac{\var^{\pi_0}( v \cdot X^0_n )}{n}
  \,.
\end{equation}

Note, when we interchanged the long time and vanishing noise limit above, we replaced the arbitrary initial distribution~$\mu_0$ on the left, with the uniform distribution~$\pi_0$.
Indeed, when~$\epsilon > 0$ the mixing effects of the noise eliminate the dependence of the asymptotic variance on the initial distribution, and the proof of Theorem~\ref{t:residual-diffusivity} in fact shows
\begin{equation}\label{e:residual-diffusivity3}
  \lim_{\epsilon \to 0}
  \lim_{n \to \infty} \frac{\var^{\mu_0}( v \cdot X^\epsilon_n )}{n}
  = \lim_{\epsilon \to 0} \lim_{n \to \infty} \frac{\var^{\pi_0}( v \cdot X^\epsilon_n )}{n}
  \,,
\end{equation}
for a sufficiently mixing map~$\varphi$, and any subgaussian initial distribution~$\mu_0$.
In general, the right-hand side of~\eqref{e:residual-diffusivity3} is not any easier to compute than the left hand side.
However, for the class of Bernoulli maps we consider, have a nice property that enables us to compute the right-hand side of~\eqref{e:residual-diffusivity3} explicitly, and this leads to the proof of Theorem~\ref{t:residual-diffusivity}.

Returning to~\eqref{e:residual-diffusivity2}, we mention that the right-hand side of~\eqref{e:residual-diffusivity2} has appeared in several works in dynamical systems
(see e.g.~\cite{Dolgopyat04, BernardinHuveneersEA15}).
In fact, the Green--Kubo formula shows that the asymptotic variance appearing on the right of~\eqref{e:residual-diffusivity2} satisfies
\begin{align}\label{eq:green-kubo}
  \lim_{n \to \infty} \frac{\var^{\pi_0}( v \cdot X^0_n )}{n}
  &\to \var(v \cdot (X^0_1 - X^0_0))
  \\
  &\quad+ 2\sum_{n=1}^\infty \cov \left( v \cdot (X^0_1 - X^0_0), v \cdot (X^0_{n+1} - X^0_n)\right)
\end{align}
as $\epsilon \to 0$. While this formula has several applications (e.g.\ proving a central limit theorem), we note that it is not explicit as it as it still involves long time dynamics.
Moreover the right-hand side may vanish, in which case there is no residual diffusivity.
There are, however, examples (including the expanding Bernoulli maps we consider) where one can check the right-hand side of~\eqref{e:residual-diffusivity2} is strictly positive, and these are all candidates where we expect to observe residual diffusivity, if the asymptotic variance can be proven independent of the initial distribution.



\subsection*{Plan of this paper}

In Section~\ref{s:bernoulli} we fix our notation convention and precisely state the assumptions under which Theorem~\ref{t:residual-diffusivity} is true.
In Section~\ref{s:proof} we explain the main idea behind the proof of Theorem~\ref{t:residual-diffusivity}, and carry out the details modulo the computation of the asymptotic variance when the initial distribution is uniform (Lemma~\ref{l:UnifID}, below).
Finally in Section~\ref{s:UnifID} we prove Lemma~\ref{l:UnifID}.

\subsection*{Acknowledgements}

The authors would like to thank
  Jon DeWitt,
  Albert Fannjiang,
  and
  Jack Xin,
for helpful discussions.

\section{Expanding Bernoulli maps}\label{s:bernoulli}

We begin by precisely describing the class of maps~$\varphi$ for which Theorem~\ref{t:residual-diffusivity} holds.
Partition $\R^d$ into unit cubes $\set{Q_k \st k \in \Z^d}$, where $Q_k = k + [0, 1)^d$.
Let $M \geq 2$ and $E_1, \dots, E_M \subseteq Q_0$ be a Borel measurable partition of $Q_0$, with
\begin{equation}
  \abs{E_1} \leq \abs{E_2} \cdots \leq \abs{E_M}
  \,.
\end{equation}
For each $i \in \set{1, \dots, M}$, let~$o_i \in \Z^d$ and~$\varphi_i \colon E_i \to Q_{o_i}$ be a Borel measurable bijection which pushes forward the normalized Lebesgue measure on~$E_i$ to the Lebesgue measure on~$Q_{o_i}$.

Given $x \in \R^d$ we let $n = \floor{x}$ denote the unique element in $\Z^d$ such that $x \in Q_n = n + Q_0$.
Define~$\varphi \colon \R^d \to \R^d$ by
\begin{equation}\label{e:PhiDef}
  \varphi(x) = n + \varphi_i(x - n)
  \quad \text{if } x - n \in E_i
  \,.
\end{equation}
Clearly the function $x \mapsto \varphi(x) - x$ is $\Z^d$ periodic, and so~$\varphi$ projects to a well-defined map~$\tilde \varphi \colon \T^d \to \T^d$ given by
\begin{equation}
  \tilde \varphi(\tilde x) = \tilde y
  \,,
  \qquad\text{where}\quad
  y = \varphi(x)
  \,.
\end{equation}
Here~$x, y \in \R^d$ and~$\tilde x, \tilde y \in \T^d$ denote the equivalence classes~$x$ and~$y$ modulo~$\Z^d$ respectively.  

We note the map $\tilde \varphi$ defined above is
conjugate to a one-sided Bernoulli shift on sequences $\{1,\dots,M\}^\mathbb{Z}$.
Moreover, since~$\tilde \varphi$ expands each set~$\tilde E_i$ to the torus~$\T^d$, we can view~$\tilde \varphi$ as an expanding Bernoulli map.
In addition to the above structure, we require a mixing assumption~$\tilde \varphi$ in order to prove Theorem~\ref{t:residual-diffusivity}.
We now state this assumption precisely.

\begin{assumption}\label{a:phiTilde}
  The map $\tilde \varphi\colon \T^d \to \T^d$, defined as above, is piecewise~$C^1$ and exponentially mixing.
  That is, there exist constants~$D < \infty$ and~$\gamma > 0$ such that for every pair of test functions~$\tilde f, \tilde g \in H^1(\T^d)$, we have
  \begin{equation}\label{e:correlationdecay}
    \abs[\Big]{
      \ip{\tilde f, \tilde g \circ \tilde \varphi^n} - \int_{\T^d} \tilde f \, dx \int_{\T^d} \tilde g \, dx
    }
    \leq D e^{-\gamma n} \norm{\tilde f}_{H^1} \norm{\tilde g}_{H^1}
    \,.
  \end{equation}
\end{assumption}

We clarify that~$\tilde E_i \subseteq \T^d$ is the projection of~$E_i \subseteq \R^d$ to the torus~$\T^d$.
We also note that~\eqref{e:PhiDef} and the fact that~$\tilde \varphi \colon \T^d \to \T^d$ is  Lebesgue measure preserving implies~$\varphi \colon \R^d \to \R^d$ is also Lebesgue measure preserving.  The shift invariance of~$\varphi$ implies the projected process~$\tilde X^\epsilon$ is a Markov process on the torus~$\T^d$.
  Equation~\eqref{e:Xn} and the fact that~$\varphi$ is Lebesgue measure preserving implies that the stationary distribution of~$\tilde X^\epsilon$ is $\tilde \pi$, the uniform measure on $\T^d$.

\begin{remark}
  The exponential mixing requirement in Assumption~\ref{a:phiTilde} can be weakened.
  In the proof, this condition is only used to ensure the right-hand side of~\eqref{e:varJ} in Proposition~\ref{p:varJ} (below) vanishes.
  Theorem~\ref{t:residual-diffusivity} will still hold provided we replace Assumption~\ref{a:phiTilde} with the assumption that
  \begin{equation}\label{e:EpTmixVanish}
    \lim_{\epsilon \to 0} \epsilon \ttmix^\epsilon = 0
    \,,
  \end{equation}
  where~$\ttmix^\epsilon$ is the~$1/2$ mixing time of~$\tilde X^\epsilon$ on $\T^d$.
  Available results (see for instance~\cite{FengIyer19,IyerLuEA24,IyerZhou23}) show even a quadratic decay in~\eqref{e:correlationdecay}  implies~\eqref{e:EpTmixVanish}.
\end{remark}

\subsection*{Notation and convention.}

We now briefly clarify several convention that will be used throughout this paper.
\begin{enumerate}[(1)]
  \item
    A tilde is used to denote projections onto the torus.
    That is if~$x \in \R^d$, then $\tilde x \in \T^d$ denotes the equivalence class~$x \pmod{\Z^d}$.
    Conversely if~$\tilde y \in \T^d$, we will implicitly use~$y \in \R^d$ to denote any representative of the equivalence class~$\tilde y$.

  \item
    If~$f \colon \R^d \to \R$ is a function then~$\tilde f \colon \T^d \to \R$ denotes its periodization
    \begin{equation}
      \tilde f(\tilde x) \defeq \sum_{k \in \Z^d} f( x + k )
      \,.
    \end{equation}
  Similarly, if~$\mu$ is a measure on~$\R^d$ then~$\tilde \mu$ is the measure on~$\T^d$ defined by~$\tilde \mu(\tilde E) = \sum_{k \in \Z^d} \mu(E + k)$, where~$E \subseteq \R^d$ is any Borel set such that~$\tilde E = \set{\tilde x \st x \in E}$.

  \item 
    The expectation operator has lower precedence than multiplication.
    That is, if~$X, Y$ are random variables then~$\E XY $ denotes the expectation of the product~$\E(XY)$, and~$\E X^2$ denotes the expectation of the square~$\E(X^2)$.

  \item
    When~$X$ is a Markov process we will use~$\E^\mu X_n$ to denote the expectation of~$X_n$ given~$X_0 \sim \mu$.
    When~$x \in \R^d$, we use~$\E^x X_n$ to denote the expectation of~$X_n$ given~$X_0 = x$.
    For random variables that are not associated to Markov processes, we will denote their expectation by~$\E$, without any sub or superscripts.
\end{enumerate}

\section{Proof of Theorem~\ref{t:residual-diffusivity}}\label{s:proof}

To prove Theorem~\ref{t:residual-diffusivity} we first need to show that the process~$X^\epsilon$ remains subgaussian, with norm controlled independently of~$\epsilon$.
Recall~\cite{Vershynin18} the \emph{subgaussian norm} of a random variable, denoted by~$\norm{\cdot}_{\psi_2}$, is defined by
\begin{equation}\label{e:psi2normDef}
  \norm{Y}_{\psi_2} \defeq \inf\set[\big]{ c > 0 \st \E e^{|Y|^2 / c^2} \leq 2 }
  \,.
\end{equation}
\begin{lemma}\label{l:subGaussianNorm}
  There exists a constant~$\Lambda > 0$ such that for all~$n \in \N$, and all~$\epsilon \in [0, 1]$ we have
  \begin{equation}\label{e:XnConcentration}
    \norm{X^\epsilon_n}_{\psi_2}
      \leq \frac{\Lambda}{2} \paren[\big]{ \norm{X_0}_{\psi_2}  + n }
      \,.
  \end{equation}
\end{lemma}
Note that even though Lemma~\ref{l:subGaussianNorm} gives an upper bound on~$\norm{X^\epsilon_n}_{\psi_2}$ that is $\epsilon$-independent, it is more crude than Theorem~\ref{t:residual-diffusivity}.
Indeed, the bound~\eqref{e:XnConcentration} implies the quadratic upper bound~$\var(v \cdot X^\epsilon_n) \leq C |v|^2 n^2$, whereas the conclusion of Theorem~\ref{t:residual-diffusivity} involves bounds (both upper and lower) that are linear in~$n$.

The next Lemma is a special case of Theorem~\ref{t:residual-diffusivity} when the initial distribution is uniform on~$Q_0$.
\begin{lemma}\label{l:UnifID}
  If~$\pi_0 = \unif(Q_0)$ denotes the uniform distribution on~$Q_0$ then
  \begin{equation}
    \abs[\Big]{
      \lim_{n \to \infty} \frac{\var^{\pi_0}( v \cdot X^\epsilon_n )}{n}
      - \var^{\pi_0}\paren[\big]{ v \cdot \floor{X^0_1} }
    }
      \leq  C |v|^2 \sqrt{ \epsilon \abs{\ln \epsilon}^3 }
      \,.
  \end{equation}
\end{lemma}

Lemma~\ref{l:UnifID} is the only place in the proof of Theorem~\ref{t:residual-diffusivity} which relies on the Bernoulli structure of~$\varphi$.
Even if~$\tilde \varphi$ is not an expanding Bernoulli map then one can still show that the asymptotic variance of~$v \cdot X^\epsilon$ starting from any (subgaussian) initial distribution equals the asymptotic variance of~$v \cdot X^\epsilon$ starting form~$\pi_0$.
Computing this for~$\epsilon > 0$ (or even for~$\epsilon = 0$), however, is not easy in general.
In our situation we compute it because Assumption~\ref{a:phiTilde} makes the increments~$\floor{\varphi(X^\epsilon_n)} - \floor{X^\epsilon_n}$ independent (see Lemma~\ref{l:varS}, below).

Momentarily postponing the proof of Lemmas~\ref{l:subGaussianNorm} and~\ref{l:UnifID} we will now prove Theorem~\ref{t:residual-diffusivity}.

\begin{proof}[Proof of Theorem~\ref{t:residual-diffusivity}]
  Before delving into the technical details of the proof, we will first explain the main idea.
  Fix~$\epsilon \in (0, 1)$, and let~$f_n$ denote the density of~$X^\epsilon_n$ on $\R^d$.
  We will begin by choosing
  \begin{equation}\label{e:mDef}
    m = \floor{n^{1/4}}
    \,,
  \end{equation}
  and finding a probability density function $g_m$ such that
  \begin{equation}\label{e:GmProps}
    \tilde g_m = 1
    \qquad\text{and}\qquad
    \norm{f_m - g_m}_{L^1(\R^d)} = \norm{\tilde f_m - 1}_{L^1(\T^d)}
    \,.
  \end{equation}
  Once we have~$g_m$, we define two new (coupled) Markov processes~$Y_k, Y_k'$, for $k \in \{m,\dots,n\}$, using the same evolution rule~\eqref{e:Xn}, the same noise, but different initial distributions.
  We specify the initial distributions for~$Y$ and~$Y'$ at time~$m$ by
  \begin{equation}\label{e:YmID}
    Y_m \sim g_m
    \quad\text{and}\quad
    Y_m' \sim \pi_0
    \,,
  \end{equation}
  respectively.
  Since~$\tilde g_m = 1$,~$Y_m$ and~$Y'_m$ differ by an element of the integer lattice~$\Z^d$, and the periodic structure of~$\varphi$ will preserve this difference at all later times.
  This combined with Lemma~\ref{l:UnifID} will allow us to show
  \begin{equation}\label{e:YnMinusYnPrime1}
    \abs[\big]{ \var( v \cdot Y_n)
      - \var( v \cdot Y'_n)}
      \leq \abs{v}^2 o(n)
      \,.
  \end{equation}

  The process $\tilde X^\epsilon$ is a Markov process on the torus~$\T^d$ with stationary distribution $\tilde \pi$, the uniform on $\T^d$.  Additionally, available results~\cite{FengIyer19, IyerLuEA24, IyerZhou23} can be used to show that~$\tilde X^\epsilon$ is mixing on $\T^d$.
  In fact, as we will see shortly, the mixing time of~$\tilde X^\epsilon$ is at most~$O(\abs{\ln \epsilon}^3)$.
  Combined with~\eqref{e:mDef} and~\eqref{e:GmProps}, this will show that~$\norm{f_m - g_m}_{L^1}$ is small when~$n$ is large.

  Finally, using subgaussianity and the fact that any Markov evolution induces a contraction on the laws, we will show
  \begin{equation}\label{e:XnMinusYn1}
    \abs{\var(v \cdot X^\epsilon_n) - \var(v \cdot Y_n)}
      \xrightarrow{n \to \infty} 0
      \,.
  \end{equation}
  Combining~\eqref{e:YnMinusYnPrime1}, \eqref{e:XnMinusYn1} and Lemma~\ref{l:UnifID} will conclude the proof of Theorem~\ref{t:residual-diffusivity}.
  \medskip

  We will now prove each of the above claims.
  Moreover, the proofs of~\eqref{e:YnMinusYnPrime1} and~\eqref{e:XnMinusYn1} will be quantitative and we will obtain an explicit rate at which the right-hand sides vanish as~$n \to \infty$.

  \restartsteps
  \step[Constructing the function~$g_m$]
  Order the elements of~$\Z^d$ as~$\set{0 = k_0, k_1, \dots}$.
  Fix~$x \in Q_0$, and define
  \begin{equation}
    \ell_0 = \ell_0(x) = \inf\set[\Big]{\ell \in \N \st \sum_{j=0}^{\ell} f( x + k_j) > 1 }\,.
  \end{equation}
  If~$\ell_0 = \infty$, then for every~$\ell \in \N$ we define
  \begin{equation}
    g_m(x + k_\ell) =
      \begin{cases}
	f_m(x + k_\ell ) & \ell > 0
	\\
f_m(x + k_0) + 1 - \sum_{\ell' \in \N} f_m(x + k_{\ell'}) & \ell = 0
    \,.
      \end{cases}
  \end{equation}
  If~$\ell_0 < \infty$, then we define
  \begin{equation}
    g_m(x + k_\ell ) =
      \begin{cases}
	f_m(x + k_\ell) & \ell < \ell_0
	\\
	1 - \sum_{\ell' < \ell_0} f_m(x + k_{\ell'}) & \ell = \ell_0
	\\
	0 & \ell > \ell_0
	\,.
      \end{cases}
  \end{equation}

  Clearly, for every $x \in \T^d$ we have
  \begin{equation}
    \tilde{g}_m(x)
      = \sum_{\ell=0}^\infty g_m(x + k_\ell)
      = \sum_{k \in \Z^d} g_m(x + k)
      = 1
      \,.
  \end{equation}
  Moreover, when~$\ell_0(x) = \infty$ we note~$g_m(x + k) \geq f_m(x + k)$ for all~$k \in \Z^d$.
  When~$\ell_0(x) < \infty$, we note~$0 \leq g_m(x + k) \leq f_m(x + k)$ for all~$k \in \Z^d$.
  Thus, in both cases, the sign of $g_m(x + k) - f_m(x + k)$ does not change with $k$ and we have
  \begin{align}
    \sum_{k \in \Z^d} |g_m(x + k) - f_m(x + k)|
      &= \abs[\Big]{ \sum_{k \in \Z^d} \paren{g_m(x + k) - f_m(x + k)} }
    \\
      &= \abs[\Big]{1 - \sum_{k \in \Z} f_m(x + k)}
      = \abs{1-\tilde{f}_m(x)}.
  \end{align}
  So integrating over $x \in \T^d$ shows $\norm{f_m - g_m}_{L^1(\R^d)} = \norm{\tilde f_m - 1}_{L^1(\T^d)}$.
  Thus~$g_m$ is a probability density function that satisfies both conditions in~\eqref{e:GmProps}, as desired.  Observe that this construction with the choice of $k_0 = 0$ guarantees that 
\[
g_m(x) \leq f_m(x), \quad \text{for all}\; x \notin Q_0.
\]

  \step[Proof of the bound~\eqref{e:YnMinusYnPrime1}]
  Once we have~$g_m$, we define two new (coupled) Markov processes~$Y_k, Y_k'$, for $k \in \{m,\dots,n\}$, as described above.  These processes $Y$ and $Y'$ use the same evolution rule~\eqref{e:Xn} as for $X$ (for times $k \in \{m,\dots,n\}$).  Moreover, $Y$ and $Y'$ are coupled, using the same noise for both processes; they differ only in the initial condition (at time~$m$). Specifically, let $Y_m \sim g_m$, and then define $Y_m'$ by
  \begin{equation}
    Y_m' = Y_m - I_m,
  \end{equation} 
  where
  \begin{equation}
    I_m = \lfloor Y_m \rfloor =  \sum_{k \in \Z^d} k \one_{\set{Y_m \in Q_k}} \,.
  \end{equation}  
In particular, $Y_m' \sim \pi_0$, since $\tilde g_m = 1$.  For $k \in \{m,\dots,n-1\}$, define $Y_{k+1} = \varphi(Y_k) + \epsilon \xi_{k+1}$ and $Y_{k+1}' = \varphi(Y_k') + \epsilon \xi_{k+1}$.  Because of this coupling and because~$I_m$ is integer valued, we have
  \begin{equation}\label{e:YnMinusYnPrime}
    Y_k = Y_k' + I_m\,, \quad \forall \;\; k \in \{m,\dots,n\}.
  \end{equation}

  Thus
  \begin{equation}\label{e:varYMinusYprime}
    \abs{\var(v \cdot Y_n) - \var(v \cdot (Y_n'))}
      \leq \var(v \cdot I_m) + 2\abs{\cov(v \cdot Y_n', v \cdot I_m)}\,.
  \end{equation}
  By Lemma~\ref{l:subGaussianNorm} we know
  \begin{equation}\label{e:varIm}
    \var(v \cdot I_m) \leq C \abs{v}^2 m^2
    \,.
  \end{equation}
  Moreover, since~$Y'_m \sim \unif(Q_0)$, Lemma~\ref{l:UnifID} implies there exists~$N_0 = N_0(\epsilon)$ such that
  \begin{equation}\label{e:varYPrime}
    \var(v \cdot Y'_n)
      \leq C \abs{v}^2 (n - m)
      \leq C \abs{v}^2 n
      \,,
  \end{equation}
  for all~$n \geq N_0$.
  Using~\eqref{e:varIm} and~\eqref{e:varYPrime} in~\eqref{e:varYMinusYprime} immediately implies
  \begin{equation}\label{e:varYMinusYprime2}
    \abs{\var(v \cdot Y_n) - \var(v \cdot (Y_n'))}
      \leq C \abs{v}^2 (m^2 + m \sqrt{n} )
      \leq C \abs{v}^2 n^{3/4}
      \,,
  \end{equation}
  for all~$n \geq N_0$.
  This proves~\eqref{e:YnMinusYnPrime1} as desired.
  \smallskip

  \step[Proof of the bound~\eqref{e:XnMinusYn1}] 
At time $m$, the density of $X_m$ is $f_m$, and the density of $Y_m$ is $g_m$. For $k \in \{m,\dots,n\}$, both processes evolve according to the same transition probabilities, and we wish to compare their variances at time $n$.  
  Let~$h_n$ be the density of~$Y_n$, and note
  \begin{multline}
    \abs{\var(v \cdot X_n) - \var(v \cdot Y_n)}
      \leq
	\int_{\R^d} (v \cdot x)^2 \abs{f_n(x) - h_n(x)} \, dx
    \\
      \label{e:VarXnVarYn1}
	+ \paren[\Big]{
	    \int_{\R^d} \abs{v \cdot x} \abs{f_n(x) - h_n(x)} \, dx
	  }
	  \paren[\Big]{ \E (\abs{v \cdot X_n} + \abs{v \cdot Y_n} )}
    \,.
  \end{multline} 
  We will bound each term on the right-hand side by splitting the integral into two parts, one where~$\abs{x} < n^{3/2}$ and the other where~$\abs{x} \geq n^{3/2}$.

  \substep[First term in~\eqref{e:VarXnVarYn1} when $x \geq n^{3/2}$]
  We first claim that there exists~$m_0$ (depending only on the dimension), such that for all~$m \geq m_0$  we have
  \begin{equation}\label{e:YmPsi2}
    \norm{Y_m}_{\psi_2} \leq \Lambda ( \norm{X_0}_{\psi_2} + m )
    \,.
  \end{equation}

  To see this, we note that the construction of~$g_m$ guarantees
  \begin{equation}
    g_m(x) \leq f_m(x)
    \quad\text{for all } x \not\in Q_0
    \,,
  \end{equation}
  and hence
  \begin{equation}
    \P( \abs{Y_m} \geq t ) \leq \P( \abs{X_m} \geq t )
    \quad\text{for all } t \geq \sqrt{d}
    \,.
  \end{equation}
  The definition of the subgaussian norm~\eqref{e:psi2normDef} implies
  \begin{equation}
    \P( \abs{X_m} \geq t) \leq 2 e^{ -t^2 / \norm{X_m}_{\psi_2}^2}
    \,,
  \end{equation}
  and hence for any~$\alpha > \norm{X_m}_{\psi_2}$ we have
  \begin{align}
    \E e^{\abs{Y_m}^2 / \alpha^2 }
      &= 1 + \int_0^\infty \frac{2 t}{\alpha^2} e^{t^2 / \alpha^2} \P( \abs{Y_m} \geq t )\, dt
    \\
      &\leq 1
	+ \int_0^{\sqrt{d}} \frac{2 t}{\alpha^2} e^{t^2 / \alpha^2} \, dt
	+ \int_{\sqrt{d}}^\infty \frac{2 t}{\alpha^2} e^{t^2 / \alpha^2} \P( \abs{X_m} \geq t )\, dt
    \\
      &\leq e^{d / \alpha^2}
	+ \int_{0}^\infty
	    \frac{4 t}{\alpha^2} \exp\paren[\Big]{ -t^2\paren[\Big]{ \frac{1}{\norm{X_m}_{\psi_2}^2 - \frac{1}{\alpha^2}}}} \, dt
    \\
      &\leq e^{d / \alpha^2} + \frac{2 \norm{X_m}_{\psi_2}^2}{\alpha^2 - \norm{X_m}_{\psi_2}^2}
    \, .
  \end{align}
  Choosing
  \begin{equation}
    \alpha = \Lambda ( \norm{X_0}_{\psi_2} + m )
  \end{equation}
  and using Lemma~\ref{l:subGaussianNorm} we see
  \begin{equation}
    \E e^{\abs{Y_m}^2 / \alpha^2} \leq 2 \,,
    \quad\text{provided}\quad
    \Lambda^2 ( \norm{X_0}_{\psi_2} + m )^2 \geq \frac{d}{\ln (4/3) }
    \,.
  \end{equation}
  This proves~\eqref{e:YmPsi2}, as claimed.

  Now, since the processes $Y_k$ and $X_k$ evolve by the same transition probability for $k \in \{m,\dots,n\}$, Lemma~\ref{l:subGaussianNorm} and~\eqref{e:YmPsi2} imply
  \begin{align}
    \norm{Y_n}_{\psi_2}
      \leq \frac{\Lambda}{2} (\norm{Y_m}_{\psi_2} + (n-m))
      \leq \frac{\Lambda^2}{2} (\norm{X_0}_{\psi_2} + n)
      \leq \Lambda^2 n
      \,,
  \end{align}
  provided~$n \geq \norm{X_0}_{\psi_2}$.
  This implies
  \begin{equation}
    \P( \abs{Y_n} > t ) \leq 2e^{-t^2 / (n^2 \Lambda^4) }
    \quad\text{for all } t \geq 0
    \,.
  \end{equation}
  Thus for any~$R > 0$ we have
  \begin{align}
    \E \paren[\big]{ \one_{\set{\abs{Y_n} > R}} \abs{Y_n}^2 }
      &= \int_{R}^\infty 2 t \P( \abs{Y_n} > t ) \, dt
      \leq \int_{R}^\infty 4 t e^{ - t^2 / (n^2 \Lambda^4) } \, dt
    \\
      \label{e:YnTail}
      &= 2 n^2 \Lambda^4 e^{ - R^2 / (n^2 \Lambda^4) }
      \,.
  \end{align}

  By the same argument and Lemma~\ref{l:subGaussianNorm}, we also obtain
  \begin{equation}\label{e:XnTail}
    \E \paren[\big]{ \one_{\set{\abs{X^\epsilon_n} > R}} \abs{X^\epsilon_n}^2 }
      \leq 2 n^2 \Lambda^2 e^{ - R^2 / (n^2 \Lambda^2) }
      \leq 2 n^2 \Lambda^4 e^{ - R^2 / (n^2 \Lambda^4) }
      \,.
  \end{equation}
  Choosing~$R = n^{3/2}$ and combining~\eqref{e:YnTail} and~\eqref{e:XnTail} and shows
  \begin{equation}\label{e:fh1}
    \int_{\abs{x} \geq n^{3/2}}
      (v \cdot x)^2 \abs{f_n(x) - h_n(x)} \, dx
      \leq |v|^2 4 \Lambda^2 n^2 e^{-n / \Lambda^4}
  \end{equation}
  for all~$n \geq \norm{\mu_0}_{\psi_2}$.

  \substep[First term in~\eqref{e:VarXnVarYn1} when $x < n^{3/2}$]
  We will show this term vanishes by showing~$\norm{f_n - h_n}_{L^1}$ vanishes exponentially as~$n \to \infty$.
  Let~$\ttmix^\epsilon = \ttmix^\epsilon (1/2)$ denote the mixing time of~$\tilde X^\epsilon$.
  That is, if~$\tilde p^\epsilon_k$ is the~$k$-step transition density of~$\tilde X^\epsilon$, then~$\ttmix^\epsilon$ is the smallest~$k\in \N$ for which
  \begin{equation}
    \sup_{\tilde x \in \T^d} \frac{1}{2} \norm{\tilde p^\epsilon_k(x, \cdot) - 1}_{L^1(\T^d)} \leq \frac{1}{4}
    \,.
  \end{equation}
  Using the fact that a Markov process induces an~$L^1$ contraction on the density, and the fact that (see for instance Section~4.5 in~\cite{LevinPeres17})
  \begin{equation}
    \sup_{x \in \T^d} \frac{1}{2}\norm{\tilde p^\epsilon_m(x, \cdot) - 1}_{L^1(\T^d)}
      \leq 2^{-\floor{m / \ttmix^\epsilon}}
  \end{equation}
  we see
  \begin{align}
    \int_{\abs{x} \leq n^{3/2}} (v \cdot x)^2 \abs{f_n(x) - h_n(x)} \, dx
      &\leq C \abs{v}^2 n^3 \norm{f_n - h_n}_{L^1}
      \leq C \abs{v}^2 n^3 \norm{f_m - h_m}_{L^1}
    \\
      &= C \abs{v}^2 n^3 \norm{f_m - g_m}_{L^1}
      = C \abs{v}^2 n^3 \norm{\tilde f_m - 1}_{L^1}
    \\
      \label{e:fh2}\noeqref{e:fh2}
      &\leq C \abs{v}^2 n^3 2^{-\floor{m / \ttmix^\epsilon} }
      \,.
  \end{align}

  \substep[Second integral in~\eqref{e:VarXnVarYn1}]
  Using the same argument as in the previous two steps we see
  \begin{align}
    \abs[\Big]{
	\int_{\abs{x} \geq n^{3/2} } (x \cdot v) (f_n(x) - h_n(x)) \, dx
      }
      &\leq 4 |v| \int_{n^{3/2}}^\infty e^{-t^2 / (n^2 \Lambda^4)} \, dt
    \\
      \label{e:fh3}\noeqref{e:fh3}
      &\leq 4 |v| \Lambda^4\sqrt{n} e^{-n/\Lambda^4}
  \end{align}
  and
  \begin{equation}\label{e:fh4}\noeqref{e:fh4}
    \abs[\Big]{
	\int_{\abs{x} \leq n^{3/2} } (x \cdot v) (f_n(x) - h_n(x)) \, dx
      }
    \leq  C \abs{v} n^{3/2} 2^{ - \floor{m / \ttmix^\epsilon} }
    \,.
  \end{equation}
  Moreover, by Lemma~\ref{l:subGaussianNorm} and~\eqref{e:YmPsi2} we see
  \begin{align}
    \label{e:fh5}\noeqref{e:fh5}
    \E \abs{X^\epsilon_n}
      &\leq \norm{X^\epsilon_n}_{\psi_2}
      \leq \frac{\Lambda}{2} \paren[\big]{ \norm{X_0}_{\psi_2}  + n }
      \,,
    \\
    \label{e:fh6}\noeqref{e:fh6}
    \text{and}\quad
    \E \abs{Y_n}
      &\leq \norm{Y_n}_{\psi_2}
      \leq \Lambda \paren[\big]{ \norm{X_0}_{\psi_2}  + n }
      \,.
  \end{align}

  Using~\eqref{e:fh1}--\eqref{e:fh6} in~\eqref{e:VarXnVarYn1} and increasing~$N_0$ if necessary, and recalling~\eqref{e:mDef} we obtain~\eqref{e:XnMinusYn1}, concluding Step~\arabic{GIStep}.
  \smallskip

  \step
  Combining~\eqref{e:YnMinusYnPrime1} and~\eqref{e:XnMinusYn1} we obtain
  \begin{equation}
    \lim_{n \to \infty} \frac{\var(v \cdot X_n)}{n}
      = \lim_{n \to \infty} \frac{\var(v \cdot Y_n)}{n}
      = \lim_{n \to \infty} \frac{\var(v \cdot Y_n')}{n}
      \,,
  \end{equation}
  and using Lemma~\ref{l:UnifID} concludes the proof.
\end{proof}

It remains to prove Lemmas~\ref{l:subGaussianNorm} and~\ref{l:UnifID}.
The proof of Lemma~\ref{l:subGaussianNorm} follows quickly from Hoeffding's inequality, and we present it here.
The proof of Lemma~\ref{l:UnifID} is more involved and we present it in Section~\ref{s:UnifID}, below.
\begin{proof}[Proof of Lemma~\ref{l:subGaussianNorm}]
  Notice first
  \begin{equation}
    X^\epsilon_{n+1} - X^\epsilon_n
      = \varphi(X^\epsilon_n) - X^\epsilon_n + \epsilon \xi_{n+1}
	\leq A + \epsilon \abs{\xi_{n+1}}
    \,,
  \end{equation}
  where
  \begin{equation}
    A \defeq \sup_{x \in \R^d} \abs{\varphi(x) - x}
    \,.
  \end{equation}
  Hence
  \begin{equation}
    \abs{X^\epsilon_n} \leq \abs{X^\epsilon_0} + An + \epsilon \Xi_n\,,
    \quad\text{where}\quad
    \Xi_n \defeq \sum_1^{n} \abs{\xi_k} \,.
  \end{equation}
  Inequality~\eqref{e:XnConcentration} now follows from Hoeffding's inequality (see for instance Theorem 2.6.2 in~\cite{Vershynin18}).
\end{proof}
\section{Proof of Lemma~\ref{l:UnifID}}\label{s:UnifID}

To prove Lemma~\ref{l:UnifID} we write
\begin{equation}\label{e:XSRJ}
  X^\epsilon_n = S_n + J_n + R_n
\end{equation}
where
\begin{align}
  S_n &= S_{n-1}
    + \floor{ \varphi(X^\epsilon_{n-1}) }
    - \floor{ X^\epsilon_{n-1} }
  \,,
  \\
  J_n &= J_{n-1}
    + \floor{ X^\epsilon_n }
    - \floor{ \varphi(X^\epsilon_{n-1}) } \,,
  \\
  R_n &= X_n^\epsilon - S_n - J_n
  \,,
\end{align}
with $S_0 = 0$, $J_0 = 0$, $R_0 = X_0^\epsilon$.
We will prove Lemma~\ref{l:UnifID} by obtaining a bound on the asymptotic variances of each of the processes~$S, J$ and~$R$.
The Bernoulli structure of~$\varphi$ implies the increments of~$S$ are independent, which allows us to compute the variance of~$S$ exactly.

\begin{lemma}\label{l:varS}
  For any~$v \in \R^d$, $n \geq 0$ we have
  \begin{equation}\label{e:VarVDotSn}
    \var^{\pi_0}(v \cdot S_n)
      = n \var^{\pi_0}( v \cdot \floor{X^0_1} )
      \,.
  \end{equation}
\end{lemma}

To bound the asymptotic variance of~$J$, we rely on the fact that the mixing time of~$\tilde X^\epsilon$ is~$o(1/\epsilon)$, which only requires~$\tilde \varphi$ to be sufficiently mixing, and does not rely on the Bernoulli structure.

\begin{proposition}\label{p:varJ}
  If~$\ttmix^\epsilon = \ttmix^\epsilon (1/2)$ denotes the mixing time of~$\tilde X^\epsilon$, then for all~$v \in \R^d$ and all~$\epsilon \in (0, 1)$ we have
  \begin{equation}\label{e:varJ}
    \var^{\pi_0}(v \cdot J_n) \leq C \abs{v}^2 n \epsilon \ttmix^\epsilon \,.
  \end{equation}
\end{proposition}

Momentarily postponing the proofs of Lemma~\ref{l:varS} and Proposition~\ref{p:varJ}, we now prove Lemma~\ref{l:UnifID}.
\begin{proof}[Proof of Lemma~\ref{l:UnifID}]
  For any~$v \in \R^d$, \eqref{e:XSRJ} implies
  \begin{multline}\label{e:varX}
    \abs{\var^{\pi_0}(v \cdot X^\epsilon_n) - \var^{\pi_0}(v \cdot S_n)}
      \leq C \Bigl(
	\var^{\pi_0}(v \cdot J_n) + \var^{\pi_0}(v \cdot R_n)
      \\
	+ \std^{\pi_0}(v\cdot S_n)\std^{\pi_0}(v \cdot J_n)
	+ \std^{\pi_0}(v \cdot S_n)\std^{\pi_0}(v \cdot R_n)
      \Bigr)
      \,.
  \end{multline}
  Note
  \begin{equation}
      R_{n} - R_{n-1} = X_n - \floor{X_n} - (X_{n-1} - \floor{X_{n-1}})
  \end{equation}
  and hence
  \begin{equation}
    R_n - R_0 = X_n - \floor{X_n} - (X_{0} - \floor{X_{0}})
      \in (-1, 1)^d
      \,.
  \end{equation}
  This implies~$\var^{\pi_0}(R_n) \leq C$.
  Using this, Lemma~\ref{l:varS} and Proposition~\ref{p:varJ} in~\eqref{e:varX} implies
  \begin{equation}\label{e:avar1}
    \lim_{n \to \infty} \frac{ \abs{\var^{\pi_0}(v \cdot X^\epsilon_n) - \var^{\pi_0}(v \cdot S_n)}}{n}
      \leq C \abs{v}^2 \left( \epsilon \ttmix^\epsilon + \sqrt{\epsilon \ttmix^\epsilon}\right)
      \,.
  \end{equation}
  To finish the proof we only a mixing time estimate which shows~\eqref{e:EpTmixVanish} holds.
  This can be done quickly from existing results.
  Indeed, the exponentially mixing condition in Assumption~\ref{a:phiTilde} allows us to use Corollary 2.5 in~\cite{FengIyer19} to show that the dissipation time (aka the~$L^2$ mixing time) of~$\tilde X^\epsilon$ is at most~$O(\abs{\ln \epsilon}^2)$.
  Following this Proposition 1.3 in~\cite{IyerLuEA24} (or Proposition 2.2 in~\cite{IyerZhou23}) implies
  \begin{equation}
    \ttmix^\epsilon \leq C \abs{\ln \epsilon}^3
    \,.
  \end{equation}
  This immediately implies~\eqref{e:EpTmixVanish} and using this and Lemma~\ref{l:varS} in~\eqref{e:avar1} finishes the proof.
\end{proof}

The remainder of this section is devoted to the proof of Lemma~\ref{l:varS} and Proposition~\ref{p:varJ}

\subsection{Variance of \texorpdfstring{$S$}{S} (Lemma \ref{l:varS})}

The proof of Lemma~\ref{l:varS} relies on the Bernoulli structure of~$\varphi$ to show that the increments of~$S$ are independent.

\begin{proof}[Proof of Lemma~\ref{l:varS}]
  Define
  \begin{equation}
    D_n = S_{n+1} - S_{n}
      = \floor{\varphi(X^\epsilon_n)} - \floor{X^\epsilon_n}
      \,.
  \end{equation}
  By Assumption~\ref{a:phiTilde}, for every~$\ell \in \Z^d$, the event~$\set{D_{n} = \ell}$ can be partitioned into events of the form~$\set{\tilde X^\epsilon_n \in \tilde E_i}$.
  We will now show that
  \begin{equation}\label{e:DeltaN1EqkCond}
    \P^{\pi_0}( D_{n} = k \given \tilde X^\epsilon_{n-1} \in \tilde E_i )
      = \P^{\pi_0}( D_{n} = k)
      = \P^{\pi_0}\paren[\big]{ \floor{X^0_1} = k  }
      \,.
  \end{equation}

  We will first show the last equality in~\eqref{e:DeltaN1EqkCond}.
  For this, we recall~$\tilde \pi$ is the Lebesgue measure on~$\T^d$, which is the stationary distribution for the Markov process~$\tilde X^\epsilon$.
  Since~$X_0 \sim \pi_0$ by assumption, and~$\tilde \pi_0 = \tilde \pi$, we must have~$\tilde X^\epsilon_n \sim \tilde \pi$ for every~$n \in \N$.
  Thus, the right-hand side of~\eqref{e:DeltaN1EqkCond} is given by
  \begin{align}
    \P^{\pi_0}( D_{n} = k)
      &= \P\paren[\big]{ \floor{\varphi(U)} - \floor{U} = k
	  \given \tilde U \sim \tilde \pi }
    \\
      \label{e:PDnEqk}
      &= \P^{\pi_0}\paren[\big]{ \floor{X^0_1} = k  }
      \,,
  \end{align}
  proving the last equality in~\eqref{e:DeltaN1EqkCond} as desired.

  We now compute the left hand side of~\eqref{e:DeltaN1EqkCond} and show it equals the right-hand side of~\eqref{e:DeltaN1EqkCond}.
  Since~$\tilde X^\epsilon_n \sim \tilde \pi$, Assumption~\ref{a:phiTilde} implies that conditioned on the event~$\tilde X^\epsilon_{n-1} \in \tilde E_i$, the distribution of~$\varphi(\tilde X^\epsilon_{n-1})$ is also~$\tilde \pi$.
  Since~$\xi_{n}$ is independent of~$X^\epsilon_{n-1}$, this in turn implies that conditioned on the event~$\tilde X^\epsilon_{n-1} \in \tilde E_i$, the distribution of~$\tilde X^\epsilon_{n}$ is the uniform distribution on~$\T^d$.
  Hence
  \begin{align}
    \P^{\pi_0}( D_{n} = k \given \tilde X^\epsilon_{n-1} \in \tilde E_i )
      &= \P^{\pi_0}( D_{n} = k \given \tilde X^\epsilon_{n} \sim \tilde \pi )
      = \P^{\pi_0}\paren[\big]{ \floor{X^0_1} = k  }
      \,,
  \end{align}
  which concludes the proof of~\eqref{e:DeltaN1EqkCond}.
  \smallskip

  Now, \eqref{e:DeltaN1EqkCond}  implies the increments of~$S$ are i.i.d.\ 
  with distribution
  \begin{equation}
    \P^{\pi_0}(\Delta_n = k )
      = \P( \floor{X^0_1} = k  )
    \,.
  \end{equation}
  This immediately implies~\eqref{e:VarVDotSn}, completing the proof.
\end{proof}

\subsection{Decorrelation bound for \texorpdfstring{$J$}{J} (Proposition \ref{p:varJ})}

For any~$n \in \N$ define
\begin{equation}
  \Delta_n
    = J_{n+1} - J_n
    = \floor{X^\epsilon_{n+1}} - \floor{\varphi(X^\epsilon_n)}
  \,.
\end{equation}
We first claim all moments of~$\Delta_0$ are of order~$\epsilon$.
\begin{lemma}\label{l:p-moment-bound}
  For all~$\epsilon \leq 1$, and all~$p \in [1, \infty)$ we have
  \begin{equation}\label{e:p-moment-bound}
    \E^{\pi_0} \abs{\Delta_0}^p \leq C_p \epsilon
      \,.
  \end{equation}
\end{lemma}

Next, we claim that the increments~$\Delta_n$ decorrelate rapidly.

\begin{lemma}\label{l:covJ}
  For any~$v \in \R^d$, $m, n \in \N$ we have
  \begin{equation}\label{e:covJ}
    \abs{\cov^{\pi_0}( v \cdot \Delta_m, v \cdot \Delta_{m + n + 1} )}
      \leq
	C \epsilon \abs{v}^2
	 \sup_{\tilde x \in \T^d} \norm{\tilde p^\epsilon_n( \tilde x, \cdot) - 1}_{L^1(\T^d)}
	\,,
\end{equation}
  where as before~$\tilde p^\epsilon_n$ is the~$n$-step transition density of the process~$\tilde X^\epsilon_n$.
\end{lemma}

Momentarily postponing the proofs of Lemmas~\ref{l:p-moment-bound} and~\ref{l:covJ}, we prove Proposition~\ref{p:varJ}.

\begin{proof}[Proof of Proposition~\ref{p:varJ}]
  Note
  \begin{align}
    \var^{\pi_0}( v \cdot J_{N} )
      &= \var^{\pi_0}\paren[\Big]{
	  \sum_{n=0}^{N-1} v \cdot \Delta_n
	}
    \\
      &=
	\sum_{n=0}^{N-1} \var^{\pi_0}(v \cdot \Delta_n)
	+ 2 \sum_{m = 0}^{N-1}
	  \sum_{n=0}^{N -m - 1}
	    \cov^{\pi_0}(
	      v\cdot \Delta_{m},
	      v\cdot \Delta_{m+n+1}
	    )
    \\
    \label{e:varYtmp1}
      &\leq
	N \abs{v}^2
	\E^{\pi_0} [\abs{\Delta_0}^2 + \abs{\Delta_0}]
	\paren[\Big]{
	  1
	  + C \sum_{n=1}^\infty
	    \sup_{\tilde x \in \T^d} \norm{\tilde{p}_n( \tilde x, \cdot ) - 1 }_{L^1(\T^d)}
	}
    \,.
  \end{align}
  Now let~$T = \ttmix^\epsilon$, and observe that for every~$\tilde x \in \T^d$, $n \in \N$ and~$j \in \set{0, \dots, T-1}$, we have
  \begin{equation}\label{e:pnYtildeGeom}
    \norm{\tilde{p}_{nT + j}(\tilde x, \cdot) - 1}_{L^1(\T^d)}
      \leq \frac{1}{2^n}
      \,.
  \end{equation}
  Thus the series on the right-hand side of~\eqref{e:varYtmp1} is bounded by~$2T$.
  Combining this with Lemmas~\ref{l:p-moment-bound} and~\ref{l:covJ}, we obtain~\eqref{e:varJ} as desired.
\end{proof}

It remains to prove Lemmas~\ref{l:p-moment-bound} and~\ref{l:covJ}.
We begin with the proof of Lemma~\ref{l:covJ}.
\begin{proof}[Proof of Lemma~\ref{l:covJ}]
  We first claim for any~$n \in \N$,
  \begin{equation}\label{e:Epi0Deltan}
    \E^{\pi_0} \Delta_n = \E^{\pi_0} \Delta_0
    \,.
  \end{equation}
  Clearly, the Markov property implies
  \begin{equation}\label{e:Epi0Deltan1}
    \E^{\pi_0} \Delta_n
      = \E^{\pi_0} \E^{X^\epsilon_n} \Delta_0
      = \E^{\pi^\epsilon_n}  \Delta_0
      \,,
  \end{equation}
  where~$\pi^\epsilon_n = \dist(X^\epsilon_n)$.
  Now we note note that adding an element of~$\Z^d$ to~$X^\epsilon_0$ does not change~$\Delta_n$ for any~$n \in \N$.
  Thus, if~$\mu_0, \nu_0$ are any two probability measures such that~$\tilde \mu_0 = \tilde \nu_0$, then for any~$n \in \N$,
  \begin{equation}\label{e:mu0nu0}
    \E^{\mu_0} \Delta_n
      = \E^{\nu_0} \Delta_n
      \,.
  \end{equation}
  Since the Lebesgue measure~$\tilde \pi$ is the stationary distribution of~$\tilde X^\epsilon$, we note~$\tilde \pi^\epsilon_n = \tilde \pi$.
  Hence~\eqref{e:Epi0Deltan1} and~\eqref{e:mu0nu0} imply~\eqref{e:Epi0Deltan} as desired.

  Now using the Markov property and~\eqref{e:Epi0Deltan} we compute
  \begin{align}
    \cov^{\pi_0}( v \cdot \Delta_m, v \cdot \Delta_{m + n + 1} )
      &= \E^{\pi_0}
	v \cdot (\Delta_m - \E^{\pi_0} \Delta_0)
	\, v \cdot (\Delta_{m+n+1} - \E^{\pi_0} \Delta_0)
    \\
      &= \E^{\pi_0} \E^{X^\epsilon_m}
	v \cdot (\Delta_0 - \E^{\pi_0} \Delta_0)
	\, v \cdot (\Delta_{n+1} - \E^{\pi_0} \Delta_0)
    \\
      &= \E^{\pi^\epsilon_m} 
	v \cdot (\Delta_0 - \E^{\pi_0} \Delta_0)
	\, v \cdot (\Delta_{n+1} - \E^{\pi_0} \Delta_0)
    \\
    \label{e:cov1}
      &= \E^{\pi_0}
	v \cdot (\Delta_0 - \E^{\pi_0} \Delta_0)
	\, v \cdot (\Delta_{n+1} - \E^{\pi_0} \Delta_0)
      \,,
  \end{align}
  where the last inequality followed because~$\tilde \pi^\epsilon_m = \tilde \pi$.

  Next we define the function~$f \colon \R^d \to \R$ by
  \begin{equation}\label{e:fDef}
    f(x) \defeq
      \E^x v \cdot (\Delta_0 - \E^{\pi_0} \Delta_0)
      = \E^x (
	v \cdot \paren[\big]{
	  \floor{X^\epsilon_{1}} - \floor{\varphi(X^\epsilon_0)}
	  - \E^{\pi_0} \Delta_0
	}
      \,.
  \end{equation}
  The Markov property and~\eqref{e:cov1} imply
  \begin{equation}\label{e:cov2}
    \cov^{\pi_0}( v \cdot \Delta_m, v \cdot \Delta_{m + n + 1} )
      = \E^{\pi_0} \brak[\big]{
	  v \cdot (\Delta_0 - \E^{\pi_0} \Delta_0)
	  P^\epsilon_n f(X^\epsilon_1)
	}
      \,.
  \end{equation}
  Here~$P^\epsilon_n$ is the~$n$-step transition operator, whose action of functions is given by
  \begin{equation}
    P^\epsilon_n g(x)
      = \E^x g(X^\epsilon_n)
      = \int_{\R^d} p^\epsilon_n(x, y) g(y) \, dy
    \,,
  \end{equation}
  where~$p^\epsilon_n$ is the~$n$-step transition density of~$X^\epsilon$.

  Note that for~$\tilde x, \tilde y \in \T^d$, the~$n$-step transition density of~$\tilde X^\epsilon$ is given by
  \begin{equation}
    \tilde p^\epsilon_n(\tilde x, \tilde y) = \sum_{k \in \Z^d} p^\epsilon_n( x, y + k )
    \,.
  \end{equation}
  Thus if~$\tilde P^\epsilon_n$ is the~$n$-step transition operator of~$\tilde X^\epsilon$, then the action of~$\tilde P^\epsilon_n$ on functions~$\tilde g \colon \T^d \to \R$ is given by
  \begin{align}
    \tilde P^\epsilon_n \tilde g( \tilde x) 
      &\defeq \E^{\tilde x} \tilde g( \tilde X^\epsilon_n )
	= \int_{\T^d} \tilde p^\epsilon_n( \tilde x, \tilde y) \tilde g( \tilde y ) \, d\tilde y
      \\
      \label{e:tildePep}
	&= \int_{\R^d} p^\epsilon_n( x, y) g( y ) \, dy
	= P^\epsilon_n g(x)
	\,.
  \end{align}

  Returning to~\eqref{e:cov2}, we note that shift invariance of~$\varphi$ implies the function~$f$ is~$\Z^d$ periodic.
  Thus, defining~$\tilde f \colon \T^d \to \R$ by~$\tilde f(\tilde x) = f(x)$ and using~\eqref{e:cov2} and~\eqref{e:tildePep} we see
  \begin{align}
    \MoveEqLeft
    \cov^{\pi_0}( v \cdot \Delta_m, v \cdot \Delta_{m + n + 1} )
      = \E^{\pi_0} \brak[\big]{
	  v \cdot (\Delta_0 - \E^{\pi_0} \Delta_0)
	  \tilde P^\epsilon_n \tilde f(\tilde X^\epsilon_1)
	}
      \\
      &= \E^{\pi_0} \brak[\bigg]{
	  v \cdot (\Delta_0 - \E^{\pi_0} \Delta_0)
	  \paren[\Big]{
	    \tilde P^\epsilon_n \tilde f(\tilde X^\epsilon_1)
	    - \int_{\T^d} \tilde f \, d \tilde \pi
	  }
	}
      \\
      &= \E^{\pi_0} \brak[\bigg]{
	  v \cdot (\Delta_0 - \E^{\pi_0} \Delta_0)
	  \paren[\Big]{
	    \int_{\T^d}
	    \paren[\big]{ \tilde p^\epsilon_n( \tilde X^\epsilon_1, y ) - 1  }
	    \tilde f(y) \, d \tilde y
	  }
	}
      \,.
  \end{align}
  Hence
  \begin{align}
    \MoveEqLeft
    \abs[\big]{\cov^{\pi_0}( v \cdot \Delta_m, v \cdot \Delta_{m + n + 1} )}
  \\
      &\leq
	\E^{\pi_0} \abs[\big]{
	  v \cdot (\Delta_0 - \E^{\pi_0} \Delta_0)
	}
	\sup_{\tilde x \in \T^d}
	  \norm{ \tilde p^\epsilon_n( \tilde X^\epsilon_1, y ) - 1  }_{L^1(\T^d)}
	  \norm{\tilde f}_{L^\infty(\T^d)}
    \\
    \label{e:cov3}
      &=
	\norm{\tilde f}_{L^1(\T^d)}
	\sup_{\tilde x \in \T^d}
	  \norm{ \tilde p^\epsilon_n( \tilde X^\epsilon_1, y ) - 1  }_{L^1(\T^d)}
	  \norm{\tilde f}_{L^\infty(\T^d)}
      \,.
  \end{align}

  It remains to estimate~$\norm{\tilde f}_{L^1}$ and~$\norm{\tilde f}_{L^\infty}$.
  Using~\eqref{e:fDef} and Lemma~\ref{l:p-moment-bound} with~$p = 1$ we see
  \begin{equation}\label{e:f1}
    \norm{\tilde f}_{L^1}
      \leq 2 \abs{v} \E^{\pi_0} \abs{\Delta_0}
      \leq C \epsilon
      \,,
  \end{equation}
  and
  \begin{align}
    \norm{\tilde f}_{L^\infty}
      &\leq \abs{v} \paren[\Big]{
	\abs{\E^{\pi_0} \Delta_0}
	+ \sup_{x \in Q_0} \E \abs[\big]{\floor{\varphi(x) + \epsilon \xi_1} - \floor{\varphi(x)})}
      }
    \\
      \label{e:f2}
      &\leq C(1 + \epsilon) \abs{v}
	\leq C \abs{v}
      \,.
  \end{align}
  Using~\eqref{e:f1} and~\eqref{e:f2} in~\eqref{e:cov3} implies~\eqref{e:covJ}, concluding the proof.
\end{proof}

Finally, we prove Lemma~\ref{l:p-moment-bound}.
\begin{proof}[Proof of Lemma~\ref{l:p-moment-bound}]
  For any~$j \in \set{1, \dots, d}$ define
  \begin{equation}
    \Delta_n^j \defeq \e_j \cdot \Delta_n
    \,,
  \end{equation}
  and note
  \begin{equation}
    \E \abs{\Delta_0^j}^p
       = \int_{y \in Q_0} \E\abs{ \e_j \cdot \paren{\floor{\varphi(y) + \epsilon \xi_1} - \floor{\varphi(y)}}}^p \, dy
       \,.
  \end{equation}
  Clearly,
  \begin{equation}\label{e:floorPhi1}
    \abs{\e_j \cdot \paren{ \floor{\varphi(y) + \epsilon \xi_1} - \floor{\varphi(y)}} }
      \leq \one_{\set{ \e_j \cdot \epsilon \xi_1 \not\in [-\psi_j, 1 - \psi_j) }} + \epsilon \abs{\e_j \cdot \xi_1}
      \,,
  \end{equation}
  where
  \begin{equation}
    \psi_j = \e_j \cdot (\varphi(y) - \floor{\varphi(y)})
    \,.
  \end{equation}
  Hence,
  \begin{equation}\label{e:tmpEjPsiPlus1}
    \E \abs{\e_j \cdot \paren{ \floor{\varphi(y) + \epsilon \xi_1} - \floor{\varphi(y)}} }^p
      \leq C_p \paren[\big]{ \P\paren[\big]{ \e_j \cdot \epsilon \xi_1 \not\in [-\psi_j, 1 - \psi_j) } + \epsilon^p }
  \end{equation}

  When~$\psi_j \leq 1/2$ a standard Gaussian tail bound implies
  \begin{equation}
    \P \paren[\big]{ \e_j \cdot \epsilon \xi_1 \not\in [-\psi_j, 1 - \psi_j) }
     \leq 2 \P \paren[\big]{ \e_j \cdot \epsilon \xi_1 \leq -\psi_j }
     \leq 2 e^{-\psi_j^2 / (2\epsilon^2) }
     \,.
  \end{equation}
  Similarly, when~$\psi_j \geq 1/2$, we note
  \begin{equation}
    \P \paren[\big]{ \e_j \cdot \epsilon \xi_1 \not\in [-\psi_j, 1 - \psi_j) }
     \leq 2 \P \paren[\big]{ \e_j \cdot \epsilon \xi_1 \geq 1-\psi_j }
     \leq 2 e^{-(1-\psi_j)^2 / (2\epsilon^2) }
     \,.
  \end{equation}
  Using these two inequalities and integrating~\eqref{e:tmpEjPsiPlus1} in~$y$ immediately implies~\eqref{e:p-moment-bound} as desired.
\end{proof}

\bibliographystyle{halpha-abbrv}
\bibliography{gautam-refs1,gautam-refs2,preprints}

\newcommand{\etalchar}[1]{$^{#1}$}
\begin{thebibliography}{MHSW22}
\expandafter\ifx\csname url\endcsname\relax
  \def\url#1{\texttt{#1}}\fi
\expandafter\ifx\csname doi\endcsname\relax
  \def\doi#1{\burlalt{doi:#1}{http://dx.doi.org/#1}}\fi
\expandafter\ifx\csname urlprefix\endcsname\relax\def\urlprefix{URL }\fi
\expandafter\ifx\csname href\endcsname\relax
  \def\href#1#2{#2}\fi
\expandafter\ifx\csname burlalt\endcsname\relax
  \def\burlalt#1#2{\href{#2}{#1}}\fi

\bibitem[Bak11]{Bakhtin11}
Y.~Bakhtin.
\newblock Noisy heteroclinic networks.
\newblock {\em Probab. Theory Related Fields}, 150(1-2):1--42, 2011.
\newblock \doi{10.1007/s00440-010-0264-0}.

\bibitem[BBPS19]{BedrossianBlumenthalEA19}
J.~Bedrossian, A.~Blumenthal, and S.~Punshon-Smith.
\newblock Almost-sure exponential mixing of passive scalars by the stochastic
  navier-stokes equations, 2019,
  \burlalt{1905.03869}{http://arxiv.org/abs/1905.03869}.

\bibitem[BCVV95]{BiferaleCrisantiEA95}
L.~Biferale, A.~Crisanti, M.~Vergassola, and A.~Vulpiani.
\newblock Eddy diffusivities in scalar transport.
\newblock {\em Physics of Fluids}, 7(11):2725--2734, 1995.
\newblock \doi{10.1063/1.868651}.

\bibitem[BCZG23]{BlumenthalCotiZelatiEA23}
A.~Blumenthal, M.~Coti~Zelati, and R.~S. Gvalani.
\newblock Exponential mixing for random dynamical systems and an example of
  {P}ierrehumbert.
\newblock {\em Ann. Probab.}, 51(4):1559--1601, 2023.
\newblock \doi{10.1214/23-aop1627}.

\bibitem[BHL{\etalchar{+}}15]{BernardinHuveneersEA15}
C.~Bernardin, F.~c. Huveneers, J.~L. Lebowitz, C.~Liverani, and S.~Olla.
\newblock Green-{K}ubo formula for weakly coupled systems with noise.
\newblock {\em Comm. Math. Phys.}, 334(3):1377--1412, 2015.
\newblock \doi{10.1007/s00220-014-2206-7}.

\bibitem[BLP78]{BensoussanLionsEA78}
A.~Bensoussan, J.-L. Lions, and G.~Papanicolaou.
\newblock {\em Asymptotic analysis for periodic structures}, volume~5 of {\em
  Studies in Mathematics and its Applications}.
\newblock North-Holland Publishing Co., Amsterdam-New York, 1978.

\bibitem[CFIN23]{ChristieFengEA23}
A.~Christie, Y.~Feng, G.~Iyer, and A.~Novikov.
\newblock Speeding up {L}angevin dynamics by mixing, 2023,
  \burlalt{2303.18168}{http://arxiv.org/abs/2303.18168}.

\bibitem[CS89]{ChildressSoward89}
S.~Childress and A.~M. Soward.
\newblock Scalar transport and alpha-effect for a family of cat's-eye flows.
\newblock {\em J. Fluid Mech.}, 205:99--133, 1989.
\newblock \doi{10.1017/S0022112089001965}.

\bibitem[Dol98]{Dolgopyat98}
D.~Dolgopyat.
\newblock On decay of correlations in {A}nosov flows.
\newblock {\em Ann. of Math. (2)}, 147(2):357--390, 1998.
\newblock \doi{10.2307/121012}.

\bibitem[Dol04]{Dolgopyat04}
D.~Dolgopyat.
\newblock Limit theorems for partially hyperbolic systems.
\newblock {\em Trans. Amer. Math. Soc.}, 356(4):1637--1689, 2004.
\newblock \doi{10.1090/S0002-9947-03-03335-X}.

\bibitem[ELM23]{ElgindiLissEA23}
T.~M. Elgindi, K.~Liss, and J.~C. Mattingly.
\newblock Optimal enhanced dissipation and mixing for a time-periodic,
  lipschitz velocity field on $\mathbb{T}^2$, 2023,
  \burlalt{2304.05374}{http://arxiv.org/abs/2304.05374}.

\bibitem[EZ19]{ElgindiZlatos19}
T.~M. Elgindi and A.~Zlato\v{s}.
\newblock Universal mixers in all dimensions.
\newblock {\em Adv. Math.}, 356:106807, 33, 2019.
\newblock \doi{10.1016/j.aim.2019.106807}.

\bibitem[FI19]{FengIyer19}
Y.~Feng and G.~Iyer.
\newblock Dissipation enhancement by mixing.
\newblock {\em Nonlinearity}, 32(5):1810--1851, 2019.
\newblock \doi{10.1088/1361-6544/ab0e56}.

\bibitem[FNW04]{FannjiangNonnenmacherEA04}
A.~Fannjiang, S.~Nonnenmacher, and L.~{Wo\l owski}.
\newblock Dissipation time and decay of correlations.
\newblock {\em Nonlinearity}, 17(4):1481--1508, 2004.
\newblock \doi{10.1088/0951-7715/17/4/018}.

\bibitem[FP94]{FannjiangPapanicolaou94}
A.~Fannjiang and G.~Papanicolaou.
\newblock Convection enhanced diffusion for periodic flows.
\newblock {\em SIAM J. Appl. Math.}, 54(2):333--408, 1994.
\newblock \doi{10.1137/S0036139992236785}.

\bibitem[FP97]{FannjiangPapanicolaou97}
A.~Fannjiang and G.~Papanicolaou.
\newblock Convection-enhanced diffusion for random flows.
\newblock {\em J. Statist. Phys.}, 88(5-6):1033--1076, 1997.
\newblock \doi{10.1007/BF02732425}.

\bibitem[Fre64]{Freidlin64}
M.~I. Fre{\u\i}dlin.
\newblock The {D}irichlet problem for an equation with periodic coefficients
  depending on a small parameter.
\newblock {\em Teor. Verojatnost. i Primenen.}, 9:133--139, 1964.

\bibitem[FW03]{FannjiangWoowski03}
A.~Fannjiang and L.~{Wo\l owski}.
\newblock Noise induced dissipation in {L}ebesgue-measure preserving maps on
  {$d$}-dimensional torus.
\newblock {\em J. Statist. Phys.}, 113(1-2):335--378, 2003.
\newblock \doi{10.1023/A:1025787124437}.

\bibitem[FW12]{FreidlinWentzell12}
M.~I. Freidlin and A.~D. Wentzell.
\newblock {\em Random perturbations of dynamical systems}, volume 260 of {\em
  Grundlehren der Mathematischen Wissenschaften [Fundamental Principles of
  Mathematical Sciences]}.
\newblock Springer, Heidelberg, third edition, 2012.
\newblock \doi{10.1007/978-3-642-25847-3}.
\newblock Translated from the 1979 Russian original by Joseph Sz{\"u}cs.

\bibitem[Hei03]{Heinze03}
S.~Heinze.
\newblock Diffusion-advection in cellular flows with large {P}eclet numbers.
\newblock {\em Arch. Ration. Mech. Anal.}, 168(4):329--342, 2003.
\newblock \doi{10.1007/s00205-003-0256-7}.

\bibitem[HIK{\etalchar{+}}18]{HairerIyerEA18}
M.~Hairer, G.~Iyer, L.~Koralov, A.~Novikov, and Z.~Pajor-Gyulai.
\newblock A fractional kinetic process describing the intermediate time
  behaviour of cellular flows.
\newblock {\em Ann. Probab.}, 46(2):897--955, 2018.
\newblock \doi{10.1214/17-AOP1196}.

\bibitem[HKPG16]{HairerKoralovEA16}
M.~Hairer, L.~Koralov, and Z.~Pajor-Gyulai.
\newblock From averaging to homogenization in cellular flows---an exact
  description of the transition.
\newblock {\em Ann. Inst. Henri Poincar\'e Probab. Stat.}, 52(4):1592--1613,
  2016.
\newblock \doi{10.1214/15-AIHP690}.

\bibitem[ILN24]{IyerLuEA24}
G.~Iyer, E.~Lu, and J.~Nolen.
\newblock Using {B}ernoulli maps to accelerate mixing of a random walk on the
  torus.
\newblock {\em Quart. Appl. Math.}, 82(2):359--390, 2024.
\newblock \doi{10.1090/qam/1668}.

\bibitem[IZ23]{IyerZhou23}
G.~Iyer and H.~Zhou.
\newblock Quantifying the dissipation enhancement of cellular flows.
\newblock {\em SIAM J. Math. Anal.}, 55(6):6496--6516, 2023.
\newblock \doi{10.1137/22M1524576}.

\bibitem[KLX22]{KaoLiuEA22}
C.~Kao, Y.-Y. Liu, and J.~Xin.
\newblock A semi-{L}agrangian computation of front speeds of {G}-equation in
  {ABC} and {K}olmogorov flows with estimation via ballistic orbits.
\newblock {\em Multiscale Model. Simul.}, 20(1):107--117, 2022.
\newblock \doi{10.1137/20M1387699}.

\bibitem[Kor04]{Koralov04}
L.~Koralov.
\newblock Random perturbations of 2-dimensional {H}amiltonian flows.
\newblock {\em Probab. Theory Related Fields}, 129(1):37--62, 2004.
\newblock \doi{10.1007/s00440-003-0320-0}.

\bibitem[LP17]{LevinPeres17}
D.~A. Levin and Y.~Peres.
\newblock {\em Markov chains and mixing times}.
\newblock American Mathematical Society, Providence, RI, 2017.
\newblock \doi{10.1090/mbk/107}.

\bibitem[LWXZ22]{LyuWangEA22}
J.~Lyu, Z.~Wang, J.~Xin, and Z.~Zhang.
\newblock A convergent interacting particle method and computation of {KPP}
  front speeds in chaotic flows.
\newblock {\em SIAM J. Numer. Anal.}, 60(3):1136--1167, 2022.
\newblock \doi{10.1137/21M1410786}.

\bibitem[LXY17]{LyuXinEA17}
J.~Lyu, J.~Xin, and Y.~Yu.
\newblock Computing residual diffusivity by adaptive basis learning via
  spectral method.
\newblock {\em Numer. Math. Theory Methods Appl.}, 10(2):351--372, 2017.
\newblock \doi{10.4208/nmtma.2017.s08}.

\bibitem[LXY18]{LyuXinEA18}
J.~Lyu, J.~Xin, and Y.~Yu.
\newblock Residual diffusivity in elephant random walk models with stops.
\newblock {\em Commun. Math. Sci.}, 16(7):2033--2045, 2018.
\newblock \doi{10.4310/CMS.2018.v16.n7.a12}.

\bibitem[MCX{\etalchar{+}}17]{MurphyCherkaevEA17}
N.~B. Murphy, E.~Cherkaev, J.~Xin, J.~Zhu, and K.~M. Golden.
\newblock Spectral analysis and computation of effective diffusivities in
  space-time periodic incompressible flows.
\newblock {\em Ann. Math. Sci. Appl.}, 2(1):3--66, 2017.
\newblock \doi{10.4310/AMSA.2017.v2.n1.a1}.

\bibitem[MCZ{\etalchar{+}}20]{MurphyCherkaevEA20}
N.~B. Murphy, E.~Cherkaev, J.~Zhu, J.~Xin, and K.~M. Golden.
\newblock Spectral analysis and computation for homogenization of advection
  diffusion processes in steady flows.
\newblock {\em J. Math. Phys.}, 61(1):013102, 34, 2020.
\newblock \doi{10.1063/1.5127457}.

\bibitem[MHSW22]{MyersHillSturmanEA22}
J.~Myers~Hill, R.~Sturman, and M.~C. Wilson.
\newblock Exponential mixing by orthogonal non-monotonic shears.
\newblock {\em Physica D: Nonlinear Phenomena}, 434:133224, 2022.
\newblock \doi{10.1016/j.physd.2022.133224}.

\bibitem[MK99]{MajdaKramer99}
A.~J. Majda and P.~R. Kramer.
\newblock Simplified models for turbulent diffusion: {T}heory, numerical
  modelling, and physical phenomena.
\newblock {\em Phys. Rep.}, 314(4-5):237--574, 1999.
\newblock \doi{10.1016/S0370-1573(98)00083-0}.

\bibitem[NPR05]{NovikovPapanicolaouEA05}
A.~Novikov, G.~Papanicolaou, and L.~Ryzhik.
\newblock Boundary layers for cellular flows at high {P}\'eclet numbers.
\newblock {\em Comm. Pure Appl. Math.}, 58(7):867--922, 2005.
\newblock \doi{10.1002/cpa.20058}.

\bibitem[OTD21]{OakleyThiffeaultEA21}
B.~W. Oakley, J.-L. Thiffeault, and C.~R. Doering.
\newblock On mix-norms and the rate of decay of correlations.
\newblock {\em Nonlinearity}, 34(6):3762--3782, 2021.
\newblock \doi{10.1088/1361-6544/abdbbd}.

\bibitem[Pie94]{Pierrehumbert94}
R.~T. Pierrehumbert.
\newblock Tracer microstructure in the large-eddy dominated regime.
\newblock {\em Chaos, Solitons \& Fractals}, 4(6):1091--1110, 1994.
\newblock \doi{10.1016/0960-0779(94)90139-2}.

\bibitem[PS08]{PavliotisStuart08}
G.~A. Pavliotis and A.~M. Stuart.
\newblock {\em Multiscale methods -- Averaging and homogenization}, volume~53
  of {\em Texts in Applied Mathematics}.
\newblock Springer, New York, 2008.

\bibitem[RZ07]{RyzhikZlatos07}
L.~Ryzhik and A.~Zlato\v{s}.
\newblock K{PP} pulsating front speed-up by flows.
\newblock {\em Commun. Math. Sci.}, 5(3):575--593, 2007.
\newblock \urlprefix\url{http://projecteuclid.org/euclid.cms/1188405669}.

\bibitem[SOW06]{SturmanOttinoEA06}
R.~Sturman, J.~M. Ottino, and S.~Wiggins.
\newblock {\em The mathematical foundations of mixing}, volume~22 of {\em
  Cambridge Monographs on Applied and Computational Mathematics}.
\newblock Cambridge University Press, Cambridge, 2006.
\newblock \doi{10.1017/CBO9780511618116}.

\bibitem[Tay21]{Taylor21}
G.~I. Taylor.
\newblock Diffusion by {C}ontinuous {M}ovements.
\newblock {\em Proc. London Math. Soc. (2)}, 20(3):196--212, 1921.
\newblock \doi{10.1112/plms/s2-20.1.196}.

\bibitem[Tay53]{Taylor53}
G.~Taylor.
\newblock Dispersion of soluble matter in solvent flowing slowly through a
  tube.
\newblock {\em Proc. R. Soc. Lond. A}, 219(1137):186--203, 1953.
\newblock \doi{10.1098/rspa.1953.0139}.

\bibitem[TC03]{ThiffeaultChildress03}
J.-L. Thiffeault and S.~Childress.
\newblock Chaotic mixing in a torus map.
\newblock {\em Chaos}, 13(2):502--507, 2003.
\newblock \doi{10.1063/1.1568833}.

\bibitem[Ver18]{Vershynin18}
R.~Vershynin.
\newblock {\em High-dimensional probability}, volume~47 of {\em Cambridge
  Series in Statistical and Probabilistic Mathematics}.
\newblock Cambridge University Press, Cambridge, 2018.
\newblock \doi{10.1017/9781108231596}.
\newblock An introduction with applications in data science, With a foreword by
  Sara van de Geer.

\bibitem[WXZ21]{WangXinEA21}
Z.~Wang, J.~Xin, and Z.~Zhang.
\newblock Sharp error estimates on a stochastic structure-preserving scheme in
  computing effective diffusivity of 3{D} chaotic flows.
\newblock {\em Multiscale Model. Simul.}, 19(3):1167--1189, 2021.
\newblock \doi{10.1137/19M1275516}.

\bibitem[WXZ22]{WangXinEA22}
Z.~Wang, J.~Xin, and Z.~Zhang.
\newblock Computing effective diffusivities in 3{D} time-dependent chaotic
  flows with a convergent {L}agrangian numerical method.
\newblock {\em ESAIM Math. Model. Numer. Anal.}, 56(5):1521--1544, 2022.
\newblock \doi{10.1051/m2an/2022049}.

\bibitem[You88]{Young88}
W.~R. Young.
\newblock Arrested shear dispersion and other models of anomalous diffusion.
\newblock {\em J. Fluid Mech.}, 193:129--149, Aug 1988.
\newblock \doi{10.1017/S0022112088002083}.

\bibitem[YPP89]{YoungPumirEA89}
W.~Young, A.~Pumir, and Y.~Pomeau.
\newblock Anomalous diffusion of tracer in convection rolls.
\newblock {\em Phys. Fluids A}, 1(3):462--469, 1989.
\newblock \doi{10.1063/1.857415}.

\bibitem[Zas02]{Zaslavsky02}
G.~M. Zaslavsky.
\newblock Chaos, fractional kinetics, and anomalous transport.
\newblock {\em Phys. Rep.}, 371(6):461--580, 2002.
\newblock \doi{10.1016/S0370-1573(02)00331-9}.

\bibitem[ZSW93]{ZaslavskyStevensEA93}
G.~M. Zaslavsky, D.~Stevens, and H.~Weitzner.
\newblock Self-similar transport in incomplete chaos.
\newblock {\em Phys. Rev. E}, 48:1683--1694, Sep 1993.
\newblock \doi{10.1103/PhysRevE.48.1683}.

\end{thebibliography}

\end{document}